\input amstex
\documentstyle {amsppt}
\font\lasek=lasy10 \chardef\kwadrat="32 
\def\kwadracik{{\lasek\kwadrat}}
\def\koniec{\lower 2pt\hbox{\kwadracik}\smallskip\smallpagebreak}
\def\Koniec{\qquad\lower 2pt\hbox{\text{\kwadracik}}}
\magnification=\magstep1
\NoRunningHeads
\def\R{\Bbb R}
\def\su{\operatorname {supp}}
\def\t{\times }
\def\Z{\Bbb Z}
\def\N{\Bbb N}
\def\a{\alpha }
\def\b{\beta}
\def\p{\partial}
\def\o{\omega }

\def\r{\rightarrow}
\def\i{\iota }
\def\ld{,\ldots,}
\def\inf{\infty}
\def\d{\operatorname {Diff}(T^n,\Cal F)_0}
\def\du{\operatorname {Diff}_U(T^n,\Cal F)_0}
\def\dz{\operatorname {Diff}_{U'}(T^n,\Cal F)_0}
\def\dv{\operatorname {Diff}_V(T^n,\Cal F)_0}
\def\dw{\operatorname {Diff}_W(T^n,\Cal F)_0}
\def\ti{\tilde}
\def\wt{\widetilde}

\def\f{\Cal F}
\def\fk{\Cal F_k}
\def\codim{\operatorname{codim}}
\def\perm{\frak S}
\def\sgn{\operatorname{sgn}}
\overfullrule=0pt
\hsize=14truecm
\topmatter
\title
On the perfectness of nontransitive groups of diffeomorphisms
\endtitle
\author
Stefan Haller \ \ \
Tomasz Rybicki
\endauthor
\thanks
The first autor is supported by the Fonds zur F\"orderung der
Wissenschaftlichen Forsch\-ung, P-10037-PHY.
The second author is supported by KBN, grant 2 P03A 024 10.
\endthanks
\date
June, 1998
\enddate

\abstract
It is proven that the identity component of the group preserving the
leaves of a generalized foliation is perfect. This shows that a well-known
simplicity theorem on the diffeomorphism group extends to the
nontransitive case.
\endabstract

\endtopmatter
\document

\head
1. Introduction
\endhead

The  goal of this paper is to show a perfectness theorem for a
class of diffeomorphism groups connected with foliations.
Throughout by foliations we shall understand generalized foliations
in the sense of P.Stefan [12].
 Given a smooth foliation $\f$ on a manifold $M$, the
symbol $\d$ stands for the identity component of the
group of all leaf preserving $C^{\infty}$-smooth
diffeomorphisms of $(M,\f)$ with compact support.

\proclaim
{Theorem 1.1} Let $\f$ be a foliation on $M$ with no leaves of dimension 0.
Then the  group $\d$ and its universal covering
$\wt{\d}$ are
perfect.
\endproclaim

Theorem 1.1  and its proof extend a well known theorem of W.P.Thurston [15]
stating that $\operatorname {Diff}(M)_0$, the identity component of the group
of all compactly supported diffeomorphisms of class $C^{\inf}$ on
a manifold $M$, is simple.
The proof  generalizes  the  case, when $M$ is
the torus (Theorem of M.R.Herman [6]), to  any manifold by a
reasoning involving the homology theory (cf.[1]). Further modifications
and completions of [10], where  the regular case has been considered, enables
us to prove Theorem 1.1. These refinements are necessary in the
generating as well as in the general case.

Note that, in general, a nontransitive group is not simple for obvious
reasons. Note as well that the perfectness implies the simplicity in
a large class of transitive groups of homeomorphisms (cf.[4]).

Notice that $\d$ belongs to the class of so-called
 modular diffeomorphism groups, see Section 2.
A clue and difficult part of the proof consist in showing the fragmentation
and deformation properties (Sect.4). The proof of it does not
appeal to foliations and holds in the whole class.

K.Fukui in [5] has shown that the abelianization  of the group stabilizing
a point is nontrivial. This simple example indicates the necessity of
the assumption in the theorem.

It is fundamental that any "isotopically connected"
diffeomorphism group determines
uniquely a foliation (cf. Section 2). It follows that the identity component of the
leaf preserving diffeomorphism group is the largest connected group
defining the foliation in question. The problem of the perfectness of
its subgroups still determining this foliation is very difficult in
particular cases. Of course, such subgroups, even
if "classical", need not be perfect in general
(see [1], or [11] in the nontransitive case). The regular case on
the ground of symplectic geometry has been studied by one of us in [11].

Throughout all manifolds, diffeomorphisms, foliations etc. are assumed
to be of the class $C^{\inf}$. Most of reasonings are neither true
 for the class $C^r$ ($r$ finite) nor for $C^{\o}$.

\head
2. Foliations and modular groups of diffeomorphisms
\endhead

Let us recall some concepts from [12].
 A {\it foliation of class $C^r$} is a partition $\Cal F$
of $M$ into weakly imbedded submanifolds (see below), called leaves, such that
the following condition holds. If $x$
belongs to a $k$-dimensional leaf, then
there is a local chart $(U,\phi )$ with $\phi (x)=0$, and $\phi (U)=V\t W$,
where $V$ is open
in $\R^k$, and $W$ is open in $\R^{n-k}$, such that if $L\in \Cal F$ then
$\phi (L|U)\cap (V\times W)=V\times l$, where $l=\{w\in W:\phi ^{-1}(0,w)\in L\}$.

 Let $L$ be a subset of a $C^r$-manifold
$M$ endowed  with a $C^r$-differentiable structure which makes it an
immersed submanifold. Then $L$ is {\it weakly imbedded} if for any
locally connected topological space $N$ and a continuous map $f: N\rightarrow M$
satisfying $f(N)\subset L$, the map $f: N\rightarrow L$ is continuous
as well. It follows that in this case such a differentiable structure
is unique.

It has been first stated in [12] that
orbits of any set of local $C^r$-diffeomor- phisms, $1\leq r\leq \omega $,
   form a  foliation. More precisely,
 a smooth mapping  $\phi $ of an open subset of $\R \times M$ into $M$ is said
 to be a $C^r$-{\it arrow} if (1) $\phi (t,.)=\phi _t$ is a local $C^r$-diffeomorphism
 for each $t$, possibly with empty domain, (2) $\phi _0=id$ on its domain,
 and (3) $\operatorname {dom}(\phi _t)\subset \operatorname {dom} (\phi _s)$
 whenever $0\leq s<t$.

 Given an arbitrary set of arrows $A$ let $A^*$ be the
 totality of local diffeomorphisms $\psi $ such that $\psi =\phi (t,.)$
 for some $\phi \in A,\ t\in \R$. Next
 $\hat A^*$ denotes the set consisting of all local diffeomorphisms being finite compositions
 of elements from $A^*$ or $(A^*)^{-1}=\{\psi ^{-1}: \psi \in A^*\}$, and of
 the identity.
Then the orbits of $\hat A^*$  are called {\it accessible} sets of $A$.

For $x\in M$ we let $ A(x),\ \bar A(x)$  be the
vector subspaces of $T_xM$ generated by
$$
 \{\dot \phi (t,y): \phi \in A, \phi _t(y)=x\},\quad
  \{d_y\psi (v): \psi \in \hat A^*,\psi (y)=x,v\in A(y)\},$$
respectively.

\proclaim
{Theorem 2.1 [12]}
Let $A$ be an arbitrary set of $C^r$-arrows on $M$. Then:

(i) Every accessible set of $A$ admits a (unique) $C^r$-differentiable
structure of a connected weakly imbedded submanifold of $M$.

(ii) The  collection of
 accessible sets  defines a  foliation $\Cal F=\Cal F(A)$.

(iii)  $\bar A(x)$ is the tangent distribution of  $\Cal F(A)$.
\endproclaim

The following "splitting" theorem will be of use.

\proclaim
{Theorem 2.2 [3, Theorem 2.1]} Let $\f$ be a foliation on $M$ and let $x$
lie on a $k$-dimensional leaf.
There exists a chart $(U,\phi)$ such that $\phi (x)=0$, $\phi (U)=V\t W$,
where $V$ (resp. $W$) is open in $\R^k$ (resp. $\R^{n-k}$), and the foliation
$\f|U$ is sent to a foliation $V\t \f_2$, where $\f_2$ is a foliation on
$W$ with a 0-dimensional leaf at 0.
\endproclaim

Let $G(M)\subset Diff^{\infty }(M)$ be any diffeomorphism group.
By a smooth path (or isotopy) in $G(M)$ we mean any family $\{f_t\}_{t\in \R}$ with
$f_t\in G(M)$ such that the map $(t,x)\mapsto f_t(x)$ is smooth.
Next,  $G(M)_0$  denotes the subgroup of all $f\in G(M)$ such that there is
a smooth path $\{f_t\}_{t\in \R}$ with $f_t=id$ for $t\leq 0$
and $f_t=f$ for $t\geq 1$, and such that each $f_t$ stabilizes outside a
fixed compact set.  Notice that $G(M)_0$ is the connected component
of $id$ if $G(M)$ is locally arcwise connected and $M$ is compact.

Given $G(M)$ the totality of $f_t$ as above constitutes a set of arrows. This set
determines uniquely a foliation.
Likewise,  the flow of a $C^r$ vector field is an
 arrow. Therefore any set of vector fields $\Cal X(M)$ defines a foliation.

Let $G(M)\subset Diff^{\infty}(M)$.
To any smooth path $f_t$  in $G(M)_0$ one can attach a family of vector fields
               $$
               \dot f_t={df_t\over dt}(f_t^{-1}). $$
Then the time-dependent family  $\dot f_t$ is a unique smooth
path in the Lie algebra corresponding to $G(M)_0$ which satisfies
the equality
$$
            {df_t \over dt}=X_t\circ f_t \quad \hbox {with}\quad {f_0=id}.
\leqno (2.1)            $$
Conversely, given a smooth family $X_t$  of compactly supported vector
fields there exists a unique solution $f_t$ of $(2.1)$.
Specifically, $f_t$ is a flow if and only if
the corresponding $X_t=X$ is time-independent.

Although
only few diffeomorphisms are elements of some flow, the simplicity theorem
of Thurston [15] states that actually $\operatorname{Diff}(M)_0$ is generated
by elements of flows.

{\bf Definition.} A Lie algebra of vector fields is called {\it modular}
if it is a $C^{\infty}(M)$ module which is $C^0$ closed.

 A group of diffeomorphisms $G(M)$ is said
to be {\it modular} if its Lie algebra (cf. [7]) $\frak g$ is modular.
Consequently, there is a one-to-one correspondence
between isotopies in $G(M)$ and in $\frak g$ given by (2.1).

\proclaim
{Lemma 2.3}
Let $V\subseteq\Cal X_c(M)$ be a $C^0$ closed $C^\infty(M)$ submodule.
For $x\in M$ set $E_x:=\{X(x):X\in V\}$. 
Then $V=\{X\in\Cal X_c(M):X(x)\in E_x\}$.
\endproclaim
{\it Proof}.
One inclusion ($\subseteq$) is trivial, we show the other one. So let
$X\in\Cal X_c(M)$ such that $X(x)\in E_x$ for all $x\in M$ and suppose
conversely $X\notin V$. Since $V$ is $C^0$-closed there exists
$\varepsilon\in C^\infty(M;\R^+)$ such that 
$$
Y\in\Cal X_c(M): \|Y(y)-X(y)\|\leq\varepsilon(y)\quad\forall y\in M
\quad\Rightarrow\quad Y\notin V
$$
The norm is with respect to some fixed Riemannian metric on $M$.
For $x\in M$ we choose $Y_x\in V$ with $X(x)=Y_x(x)$ and a neighborhood
$U_x$ of $x$ such that $\|Y_x(y)-X(y)\|\leq\varepsilon(y)$ for all 
$y\in U_x$. Since the support of $X$ is compact we find $x_1,\dotsc,x_n$
with $U_{x_1}\cup\cdots\cup U_{x_n}\supseteq\su(X)$. Finally we choose
a partition of unity $\lambda_0,\lambda_1,\dotsc,\lambda_n$ subordinate to
$\{M\setminus\su(X),U_{x_1},\dotsc,U_{x_n}\}$ (that is 
$\su(\lambda_0)\subseteq M\setminus\su(X)$, $\su(\lambda_i)\subseteq
U_{x_i}$)
and define $Y:=\sum_{i=1}^n\lambda_iY_{x_i}\in V$. 
Using $\lambda_0X=0$ we obtain
$$
\align
\textstyle
\|Y(y)-X(y)\|&\textstyle
=\big\|\sum_{i=1}^n\lambda_i(y)(Y_{x_i}(y)-X(y))\big\|
\\&\textstyle
\leq\sum_{i=1}^n\underbrace{\lambda_i(y)\|Y_{x_i}(y)-X(y)\|}_{\leq\lambda_i(y)\varepsilon(y)}
\leq\varepsilon(y)
\endalign
$$
for all $y\in M$ and therefore $Y\notin V$, a contradiction.
\Koniec

\proclaim
{Proposition 2.4}
Suppose $G(M)$ is modular and let $\{U_i\}$
be a finite family of open balls of $M$. If
$f_t$ is an isotopy in $G(M)$ such that $\overline {\bigcup _t
\su (f_t})\subset \bigcup U_i$ then there are  isotopies
$f^j_t$ supported in $U_{i(j)}$ which satisfy $f_t=f^s_t\circ \cdots \circ
f^1_t$.
\endproclaim
{\it Proof}.
Let $f_t$ be as above and let $X_t$ be the corresponding family in
$\Cal X_G(M)$. By considering $f_{(p/m)t}f^{-1}_{(p-1/m)t}$,
$p=1,\ldots ,m$, instead of $f_t$ we may assume that $f_t$ is close to
the identity.

 First we  choose a new family of open balls, $\{V_j\}_{j=1}^s$, satisfying
$\su (f_t)\subset V_1\cup \ldots \cup V_s$ for each $t$
and which is starwise finer that $\{U_i\}$, that is
$$(\forall j)\,(\exists i)\, \operatorname { star}(V_j)\subset U_{i(j)},\quad
\hbox {where}\quad \operatorname{star}(V_j)=\bigcup _{\bar V_j\cap \bar V_k\neq
\emptyset} V_k.$$

 Let $(\lambda _j)_{j=1}^s$ be
a partition of unity subordinate to $(V_j)$, and let $Y^j_t=\lambda _jX_t$.
We set
$$
X^j_t=Y^1_t+\cdots +Y^j_t,\quad j=1,\ldots ,s,
$$
and $X^0_t=0$.
Each of the smooth families $X^j_t$ integrates to an isotopy $g^j_t$
with support in $V_1\cup \ldots \cup V_j$. We get  the partition
$$ f_t=g_t^s=f_t^s\circ \cdots \circ f_t^1,$$
where $ f_t^j=g_t^j\circ (g_t^{j-1})^{-1}$, with the required
 inclusions
$ \su (f_t^j)=\su (g_t^j\circ (g_t^{j-1})^{-1})\subset
\operatorname {star}(V_j)\subset U_{i(j)}$ which hold if $f_t$ is sufficiently
small.
\Koniec

\proclaim
{Proposition 2.5} Let $G(M)$ be modular. Then $G(M)_0$ is locally
contractible in the $C^{\inf}$-topology.
\endproclaim

Indeed, it follows by a classical argument involving vector fields.

\head
3. The  Case of $Diff(T^n,\f)$
\endhead

The concept of $\Cal L$-category was introduced in [14].
Roughly speaking, an
object in this category is a quadruple $(E,B,\Cal N,\Cal S)$, where $E$ is
a Fr\'echet space, $\Cal N=(|\quad |_i)$ is an increasing sequence of norms
defining the topology of $E$, $\Cal S=(S_t),\ t>0,$ is a one-parameter family
of "approximation" operators on $E$, and $B$ is an open subset with respect
to some norm from $\Cal N$. Let $E_i$ denote the completion of $E$ with
respect to the norm $|\quad |_i$, and let $\rho _{ji}:E_j\rightarrow E_i$
be an extension if $id_E$, $j\geq i$. Then topologically $E=\lim _{\leftarrow}(E_i,\rho _{ji})$.
An interpretation of the operators $S_t$ is the following. Each $S_t$ extends
to an $S_t:E_0\rightarrow E$ and $S_t$ approximates an element from $E_0$ by
an element from $E$. The greater is $t$ the better is an approximation.

The concept of $C^r$ (weak) morphism in the $\Cal L$-category is even more
complicated. Of course, all morphisms are continuous mappings.

By means of the $\Cal L$-category one can introduce
 the notion of $\Cal L$-manifold of class $C^r,\ 1\leq r\leq \infty$. This
is a topological space endowed with an $\Cal L$-atlas, i.e. an atlas modeled
on an $\Cal L$-object in the usual way. In particular, the concept of
tangent space of $\Cal L$-manifold at a point is well defined.
The spaces of $C^r$ mappings are  clue examples of $\Cal L$-manifolds,
and the need of a generalized smooth structure on them motivated the definition
of $\Cal L$-category.

The object of our interest will be $\Cal L$-groups, that is topological groups
such that their group products and  inverse mappings are $\Cal L$-morphisms.
Of course, the diffeomorphism groups are here the main example.
In obvious way one can define also a notion of $\Cal L$-action of an $\Cal L$-group
on an $\Cal L$-manifold.

We begin with an Implicit Function
Theorem  in the case of $\Cal L$-actions (cf.[14]). Let $G$, $H$ be $\Cal  L$-groups
of class $C^r$ ($r\geq 2$) and $\Cal M$ be an $\Cal L$-manifold. Denote by
$\alpha :G\times G\rightarrow G, \beta :H\times H\rightarrow H$ the group
products and let $\Phi :G\times \Cal M \rightarrow \Cal M, \Psi :H\times \Cal M \rightarrow \Cal M$
be $\Cal L$-actions of class $C^r$. Next, let $\Delta :G\times H\times \Cal M\rightarrow \Cal M$,
be an "action" of $G\times H$, defined by
$$ \Delta (g,h,x)=\Phi (g,\Psi (h,x))$$
for $g\in G, h\in H, x\in \Cal M$. By $d\Delta $ we denote the differential
of $\Delta$ with respect to two first variables. By the chain rule one has

     $$  d\Delta (g,h,x,\hat{g},\hat{h})=d_1\Phi(g,\Psi(h,x),\hat{g})+
     d_2\Phi(g,\Psi(h,x),d_1\Psi(h,x,\hat{h})).$$
(Here we adopt the notation $\hat{g}\in T_g(G),\ \hat{x}\in T_x(\Cal M)$ and
so on.) Let us fix $x_0\in \Cal M$. By making use of the local triviality of
the tangent bundle $T\Cal M$ one can identify $T_x(\Cal M)$ with $T=T_{x_0}(\Cal M)$
for $x$ being near $x_0$. Likewise, $T_g(G)$ is identified with $T_1=T_e(G)$,
whenever $g\in G$ is near $e$, and $T_h(H)$ is identified with $T_2=T_e(H)$,
whenever $h\in H$ is near $e$. Then by applying Implicit Function Theorem
one has the following

\proclaim
{Theorem 3.1 [14, 4.2.5]}  Suppose that there  exists
 an $\Cal L$-morphism of class $C^{\infty}\ $,
$L: \Cal U\times T\rightarrow T_1\times T_2$ ,
where $\Cal U$ is a neighborhood of $e$ in $H$, such that if $L(h,\hat{x})=(\hat{g},\hat{h})$,
then
$$  d\Delta (e,e,\Psi (h,x),\hat{g}, \hat{h})=\hat{x} .$$
Then there exists a neighborhood $\Cal V$ of $x_0$ in $\Cal M$ and a weak
$\Cal L$-morphism of class $C^{\infty }$ \ $s:\Cal V\rightarrow G\times H$ such that
$\Delta (g,h,x_0)=x$ if $s(x)=(g,h)$.
\endproclaim

Now let $T^n$ be the $n$-dimensional torus.
Let $1\leq k< n$ and let $\Cal F_k$ denote the
trivial $k$-dimensional foliation on the torus $T^n$, i.e. $\Cal F_k=\{T^k\times \{pt\}\}$.
We have the canonical inclusion $\a \in T^k\hookrightarrow R_{\a }\in Diff^{\inf}(T^n,\fk)_0$,
 where
 $$R_{\a }(z_1,\ldots ,z_n)=(e^{2\pi i\a_1}z_1,\ldots ,e^{2\pi i\a_k}z_k,z_{k+1},
\ldots ,z_n).$$
Given a foliation $\f'$ on $T^{n-k}$ we set $\f=T^k\t \f'$,
so that $\f_k$ is a subfoliation of $\f$. It is
clear that
$$\operatorname {Diff}^{\inf}(T^n,\f k)_0  \subset \d .\leqno (3.1)$$

Recall that $\alpha =(\alpha _1,\ldots ,\alpha _n)\in \R^n$ satisfies the
{\it Diophantine condition} if there are small $c>0$ and large $N$ such that
for any $(n+1)$-tuple of integers $(q_0,q_1,\ldots ,q_n)$
with $(q_1,\ldots ,q_n)\neq 0$ one has

$$  |q_0+q_1\alpha _1+\cdots +q_n \alpha _n|\ >\ c(|q_1|+\cdots +|q_n|)^{-N}. $$
Next, $\a \in T^n$ verifies this condition if so does its representant
in $\R^n$. The set of all such $\a$ is dense in $T^n$.

\proclaim
{ Theorem 3.2}  Let $(M,\f)$ be such that $\f=T^k\times \f'$ and let
$\alpha \in T^k$ be Diophantine.
 There exist a neighborhood $\Cal U$ of $R_{\alpha }$
in $\d$ and a continuous map
$s:\Cal U\rightarrow \d \times T^k    $
such that $h =R_{\b}g^{-1}R_{\alpha}g$ whenever $h \in \Cal U$, and
$s(h)=(g,\b )$. Furthermore, if $h_t,\ t\in I$, is a smooth isotopy in $\Cal U$ and
$s(h_t)=(g_t,\b _t)$ then $g_t, \b _t$ depend smoothly on $t$.
\endproclaim

{\it Proof}.
The starting observation is
 that if $\a \in T^k$ is Diophantine then so
is $\a '=(\a ,0)\in T^n$.

Let $\bar \a \in \R^k$ be a representant of $\a \in T^k$. The
 components of $\bar \a$ can be chosen linearly
independent over $\Bbb Q$. It follows that $\a $ generates a dense subgroup
of $T^k$.

Let $G=T^k, H=\d$ endowed with the structure opposite to the usual.
Define actions of $G$ and $H$, respectively, on $H$ by
$$\eqalign{    \Phi (\lambda ,h)&=R_{\lambda }h,\cr
    \Psi (g,h)&=g^{-1}hg.\cr} $$
Let $\Delta:G\times H\times H\rightarrow H$ be the composition of these actions
$$   \Delta (\lambda ,g,h)=\Phi (\lambda ,\Psi (g,h))=R_{\lambda }g^{-1}hg.$$
We  use  Theorem 3.1. We have
$$   d\Delta (e,e,h,\hat{\lambda},\hat{g})=\hat{\lambda}+dh\cdot \hat{g}-\hat{g}\cdot \Delta,$$
where $\hat{\lambda}\in \R^k=T_e(T^k), \ \hat{g}\in =T_{id}(\d)$.
Consider the equation
$$  \hat{\lambda}+d(g^{-1}R_{\alpha}g)\cdot \hat{g}-\hat{g}\cdot (g^{-1}R_{\alpha}g)\ =\ \hat{h}.
$$
In view of Theorem 3.1 we have to solve this equality for given $g,h\in \d$,
$\hat h\in T_{id}(\d)$, and
with respect to the unknowns $\hat{\lambda},\ \hat{g}$. Set $\hat{f}=dg\cdot
\hat{g} \cdot g^{-1} \in \d$. Since $dR_ \alpha =id$,
we get
$$   \hat{f}-\hat{f}\cdot R_{\alpha}=dg\cdot (\hat{h}-\hat{\lambda})\cdot g^{-1}
 \leqno (3.2)$$

If $m$ be the normalized Haar measure on $T^n$, we have
$$
      \int _{T^n}dg \cdot(\hat{h}-\hat{\lambda})\cdot g^{-1}\ dm\ =\ 0.
      \leqno (3.3)$$
The equality (3.3) determines uniquely $\hat{\lambda }\in \R^k$, provided
$g $ is sufficiently near $id$ in $\d$.

Now by using the condition on $\a$ and expanding in Fourier series the
both sides of (3.2) we get as in [14] or [6]
the existence of $\hat{f}\in C^{\infty }(T^n,\R^n)$ satisfying
(3.2).  Observe that in case of class $C^r$, $r$ finite,  one could not avoid the "loss
of smoothness",i.e. $\hat{f}$ is of class $C^{r-\beta}$, $\beta$ depending on
$\alpha $, cf.[8].

Since $g$ is leaf preserving
it remains to show that $\hat{f}\in T_{id}(\d)$, i.e. $\hat{f}$
is a vector field tangent to $\f$.

First fix $p'\in T^k$ and denote $C(\cdot):=\hat{f}(p',\cdot)$.
Then fix $p=(p',p'')\in T^n$ and let $L_p$ (resp. $L^k_p$) be the leaf of $\f$ (resp. $\f_k$)
passing through $p$. Choose $(x,\bar{x},y)=(x_1\ld x_k,\bar{x}_{k+1}\ld \bar {x}_{l(p)},
y_1\ld y_{n-l(p)})$, a distinguished chart  of $\f$ at $p$, where
$l(p)=dim(L_p)\geq k$. Consequently, $L_p$ is determined by $y=0$ in
this chart. In addition, we may and do assume that
$\bar{x}=y=0$ describes $L_p^k$ in the domain of $(x,\bar{x},y)$.

The equation (3.2) can be rewritten as
$$
\hat{f}_i-\hat{f}_i\cdot R_{\a}=\hat{k}_i,\quad i=1\ld n. \leqno (3.4)
$$
Here $\hat{k}=(\hat{k}_i)\in T_{id}(\d)$ denotes the r.h.s. of (3.2). Observe that
$\hat{k}_i(q)=0$ if $i>dim(L_p)$, where $q\in L_p$.

W introduce $\bar{f}=(\bar{f}_i)$ by
$$\bar{f}:=\hat{f}-C,
$$
where $C$ is as above.
 Hence $\bar{f}(0,\bar{x},y)=0,\ \forall \bar{x},y$, and
$(3.4)$ is satisfied with $\bar{f}_i$ instead of $\hat{f}_i$.

For $i>k(p)$   we  have by (3.4)
 $\bar{f}_i(\a,\bar{x},0)=0$ in the chart $(x,\bar{x},y)$ , and inductively
$$
\bar{f}_i(m\a,\bar{x},0)=0,\,\, \forall \bar{x}\,\forall m\in \Z^k.$$
Therefore, by the assumption on $\a$, $\bar{f}_i|_{L_p}=0$.
By repeating this argument to any leaf of $\f$ we get
that $\bar{f}$ is tangent to $\f$. Since $C$,
and consequently $\bar{f}$, are $C^{\infty}$, the required solution
of (3.2) is $\bar{f}$.

The second assertion follows from the fact that $s$ is a
$C^{\infty }$ $\Cal L$-morphism (Theorem 3.1), and it sends smooth curves to smooth
curves, cf.[1].                     \koniec

\head
4. Topological background
\endhead

Let us fix notation and recall some facts from the homology of groups (see
e.g. [2]).
Let $G$ be a connected topological group. We shall be concerned with $H_1(G)$
which is identified with the abelianization $G/[G,G]$.

By $\tilde G$ we denote the universal covering of $G$.
Provided $G$ is locally arcwise connected,  $\tilde G$
can be thought of as the set of pairs $(g,\{g_t\})$, where $g\in G,\, g_t\, (t\in I)$
is a path connecting $g$ with $e$, and $\{g_t\}$ is the homotopy class of
$g_t$ rel. endpoints. Then $\tilde G$ is given a group structure by the
pointwise multiplication or, equivalently, by the juxtaposition.
Clearly the perfectness of $\ti G$ yields the perfectness of $G$.

With any $G$ we can associate some simplicial set  $B\bar G=\bigcup B_m\bar G$
where $B_m\bar G$ is identified with
the set $(G,e)^{(\Delta ^m,e_0)}$ of continuous mappings of the
standard $m$-simplex $\Delta ^m$ into
$G$ sending $e_0$ (the first vertex) to $e$.
For the detailed definition see [2] or [9].

We recall some properties of
 $B\bar G$. It is a Kan complex and one can give
a purely combinatorial definition of homotopy groups (cf.[9]). Namely the
following equivalence relation is given on $B\bar G$: for any 1-simplices $\sigma ,\tau
\in B_1\bar G$
$$   \sigma \sim \tau \quad \hbox {iff} \quad \exists c\in B_2\bar G:
\partial _0 c=\sigma, \partial _1 c=\tau , \partial _2 c=e  $$
where $e$ is the constant map. Then the first homotopy group  of $B\bar G$ is defined
by $\pi_1(B\bar G)=B_1\bar G/\sim $.
It follows  that for any $\sigma \in B_1\bar G$ the classes of $\sigma $
with respect to the relation $\sim $ and with respect to the homotopy rel.
endpoints are the same, that is
$   \pi _1(B\bar G)\ =\ B_1\bar G/\sim\ =\ \tilde G$.
One then has
$$   H_1(B\bar G,\Z)\ =\ H_1(\pi _1(B\bar G))\ =\ H_1(\tilde G)
$$
since $H_1(B\bar G,\Z)\ =\ \pi_1(B\bar G)/[\pi_1(B\bar G),\pi_1(B\bar G)]$.

\head
5. Fragmentation and Deformation for Modular Groups of Diffeomorphisms
\endhead

Let $G=G(M)$ be a modular group and $\frak g$ its Lie algebra.

Let $\Delta^p:=\{(t_0,\dotsc,t_p):0\leq t_i\leq
1,\sum_{i=0}^pt_i=1\}\subseteq\R^{p+1}$
denote the standard $p$-simplex and $e_0,\dotsc,e_p$ the unit vectors in
$\R^{p+1}$, that is the edges of $\Delta^p$.
If $g:\Delta^p\to G$ is smooth then the right
logarithmic derivative $\delta^rg\in\Omega^1(\Delta^p;\frak g)$
satisfies the Maurer Cartan equation 
$$
\textstyle
d(\delta^rg)-\frac{1}{2}[\delta^rg,\delta^rg]=0
$$
Conversely, if $\sigma\in\Omega^1(\Delta^p;\frak g)$ satisfies
$d\sigma-\frac{1}{2}[\sigma,\sigma]$ then there exists
$g:\Delta^p\to G$ with $\delta^rg=\sigma$.
Moreover $g$ is unique if one assumes $g(e_0)=\operatorname{id}$. (See
\cite{7}).

If $\sigma\in\Omega^1(\Delta^p;\frak g)$ then one obtains a distribution on
$\Delta^p\times M$ by setting $E_{(t,x)}=\{(X_t,\sigma(X_t)(x)):X_t\in
T_t\Delta^p\}\subseteq T_{(t,x)}\Delta^p\times M$. This distribution has 
codimension $\dim(M)$ and it is transversal
to $\{pt\}\times M$. Moreover it is integrable if and only if
$d\sigma-\frac{1}{2}[\sigma,\sigma]=0$.

So we obtain a one-to-one correspondence between
$C^\infty\big((\Delta^p,e_0),(G,e)\big)$, one forms
$\sigma\in\Omega^1(\Delta^p;\frak g)$ satisfying
$d\sigma-\frac{1}{2}[\sigma,\sigma]=0$, and foliations on $\Delta^p\times M$
of codimension $\dim(M)$ that are transversal to $\{pt\}\times M$.
We will denote the set of all such simplices by $S_p(B\bar G)$. If $\Cal U$
is a set of open sets in $M$ then $S_p^{\Cal U}(B\bar G)$ will denote the
set of all simplices with support in one of the sets of $\Cal U$. On the
free Abelian group $C_*(B\bar G;\Z)$ generated by $S_*(B\bar G)$ we have a
differential $\partial=\sum_{i=0}^p(-1)^i\partial_i$, where
$\partial_i=\delta_i^*$ and $\delta_i:\Delta^p\to\Delta^{p+1}$ are the
inclusion of the faces. Notice that $0=\partial:C_1(B\bar G;\Z)\to 
C_0(B\bar G;\Z)\cong\Z$.

For an open neighborhood $\Cal E$ of $0\in\frak g$ we let
$$
S^{\Cal E}_p(B\bar G)=
\{\sigma\in S_p(B\bar G):\sigma(D\Delta^p)\subseteq{\Cal E}\}
$$
where $D\Delta^p\subseteq T\Delta^p$ is the unit disk bundle of $\Delta^p$.
Then $C^{\Cal E}_*(B\bar G;\Z)$ is a subcomplex of $C_*(B\bar G;\Z)$.
Moreover we will consider $S^{\Cal E,\Cal U}_*:=S^{\Cal E}_*(B\bar G)
\cap S^{\Cal U}_*(B\bar G)$ and $C^{\Cal E,\Cal U}_*(B\bar G;\Z)$.
Using barycentric subdivision one can easily show 

\proclaim
{Lemma 5.1}
For any set $\Cal U$ of open sets in $M$ and for every neighborhood $\Cal E$
of $0\in\frak g$ the inclusion induces an isomorphism
$H^{\Cal E,\Cal U}_*(B\bar G;\Z)\cong H^{\Cal U}_*(B\bar G;\Z)$.
\endproclaim

\proclaim
{Lemma 5.2}
Let $\tau:\Delta^p\times M\to\Delta^q\times M$ be smooth with
$pr_M\circ\tau=pr_M$ and let $G$ be modular. If $\sigma\in
S_q(B\bar G)$ such that the foliation corresponding to $\sigma$ is
transversal to $\tau(t,\cdot):M\to\Delta^q\times M$ for all $t\in\Delta^p$ 
then $\tau^*\sigma\in S_p(B\bar G)$. Moreover we have
$\su(\tau^*\sigma)\subseteq\su(\sigma)$.
\endproclaim

{\it Proof}.
Obviously $\tau^*\sigma$ is a foliation on $\Delta^p\times M$ with
$\codim(\tau^*\sigma)=\dim(M)$ which is
transversal to $\{pt\}\times M$ and so we obtain at least 
$\tau^*\sigma\in S_p(B\overline{\operatorname{Diff}}^\infty_c(M)_\circ)$. If
$Y\in T_t\Delta^p$ the defining equation for $\tau^*\sigma(Y)$ is
$$
\align
\sigma\big(T_{(t,x)}(pr_{\Delta^q}\circ\tau)\cdot(Y,\tau^*\sigma(Y)(x))\big)(x)
&=T_{(t,x)}(pr_M\circ\tau)\cdot(Y,\tau^*\sigma(Y)(x))
\\&
=(\tau^*\sigma)(Y)(x)
\endalign
$$
So we see that $(\tau^*\sigma)(Y)(x)\in E_x:=\{X(x):X\in\frak g\}$ for all
$x\in M$ hence by Lemma~2.3 we obtain $(\tau^*\sigma)(Y)\in\frak g$ and
thus $\tau^*\sigma\in S_p(B\bar G)$.
\koniec

\proclaim
{Lemma 5.3}
Let $\tau:M\to\Delta^p$ be smooth and define
$$
{\Cal E}_\tau:=\{X\in\frak g:\|T_x\tau\cdot X_x\|<1\quad\forall x\in M\}
\subseteq\frak g.
$$ 
Then ${\Cal E}_\tau$ is a zero neighborhood in $\frak g$ and for
$\sigma\in S_p^{{\Cal E}_\tau}(B\bar G)$ the foliation on 
$\Delta^p\times M$ corresponding to $\sigma$ is transversal to
$(\tau,\operatorname{id}_M):M\to\Delta^p\times M$.
\endproclaim

\proclaim
{Corollary 5.4}
Let $\tau_i:M\to\Delta^p$ for $i=1,\dotsc,N$. Then there exists a zero
neighborhood ${\Cal E}$ in $\frak g$ such that for
$\sigma\in S^{\Cal E}_p(B\bar G)$ the
foliation on $\Delta^p\times M$ corresponding to $\sigma$ is transversal 
to $(\mu,\operatorname{id}_M)$, where $\mu:=\sum_{i=1}^Nt_i\tau_i$ 
is any convex combination of the $\tau_i$, that is $0\leq t_i\leq 1$ and
$\sum_{i=1}^Nt_i=1$.
\endproclaim

For $N\in\N$ let
$D^n_N:=\{(m_0,\dotsc,m_n)\in\N_0^{n+1}:\sum_{i=0}^nm_i=N\}$.
If $\lambda\in\Delta^{N-1}$ and $m\in D^n_N$ we define
$\tau^\lambda_m\in\Delta^n$ by
$$
\tau^\lambda_m:=(\underbrace{\lambda_0+\cdots+\lambda_{m_0-1}}_{m_0},
\underbrace{\lambda_{m_0}+\cdots+\lambda_{m_0+m_1-1}}_{m_1},\dotsc,
\underbrace{\dotsc}_{m_n})
$$
Moreover we let
$$
A^n_N:=\{(m,\pi)\in D^n_N\times\perm_n:m+f_{\pi(1)}+\cdots+f_{\pi(j)}\in D^n_N\quad\forall 0\leq j\leq
n\}
$$
where $f_i:=e_i-e_{i-1}$. If $(m,\pi)\in A^n_N$ we define
$\tau^\lambda_{(m,\pi)}:\Delta^n\to\Delta^n$ by
$\tau^\lambda_{(m,\pi)}(e_j):=\tau^\lambda_{m+f_{\pi(1)}+\dotsc+f_{\pi(j)}}$
for $0\leq j\leq n$ and extend it affine.

If $\lambda\in C^\infty(M,\Delta^{N-1})$ is a finite partition of unity and
$(m,\pi)\in A^n_N)$ we define 
$$
\tau^\lambda_{(m,\pi)}:\Delta^n\times M\to\Delta^n\times M
\qquad\tau^\lambda_{(m,\pi)}(t,x)=(\tau^{\lambda(x)}_{(m,\pi)}(t),x)
$$
>From Corollary~5.4 we obtain a zero neighborhood $\Cal E^\lambda_n$ in
such that $\tau^\lambda_{(m,\pi)}(t,\cdot):M\to\Delta^p\times M$ is
transversal to the foliation corresponding to $\sigma\in S^{\Cal E^\lambda_n}_p(B\bar
G)$ for all $p\leq n$ and $(m,\pi)\in A^p_N$. Moreover Lemma~5.2 yields
$$
\sigma\in S^{\Cal E^\lambda_n,\Cal U}_p(B\bar G)
\Rightarrow
(\tau^\lambda_{(m,\pi)})^*\sigma\in S_p^{\Cal U}(B\bar G)
$$
for all $p\leq n$ and $(m,\pi)\in A^p_N$. We are now in a position to give
the following

\definition
{Definition 5.5 (Fragmentation mapping)}
Let $G$ be a modular,
$\Cal U$ a set of open sets in $M$, $N\in\N$ and 
$\lambda\in C^\infty(M,\Delta^{N-1})$. Then for $p\leq n$ we define 
$$
\varphi_p^\lambda:C_p^{{\Cal E}_n^\lambda,{\Cal U}}(B\bar G;\Z) 
\to C_p^{\Cal U}(B\bar G;\Z)
\qquad
\varphi^\lambda_p(\sigma):=\sum\limits_{(m,\pi)\in
A^p_N}\sgn(\pi)(\tau^\lambda_{(m,\pi)})^*\sigma
$$
where the simplex 
$\sigma$ is considered as foliation on $\Delta^p\times M$.
\enddefinition

This is a modification of the A.~Banyaga's procedure  used for the
deformation for globally hamiltonian diffeomorphisms, see \cite{1}.

Next we subdivide $\Delta^p\times I$ into $p+1$ simplexes in the
usual way. For $1\leq i\leq p+1$ we define 
$s_i^{p+1}:\Delta^{p+1}\to\Delta^p\times I$ by
$$
s_i^{p+1}(e_j):=
\cases
(e_j,0) & 0\leq j<i
\\
(e_{j-1},1) & i\leq j\leq p+1
\endcases
$$
and extend it affine. 
For $(m,\pi)\in A^p_N$ we define $T_{(m,\pi)}^\lambda:\Delta^p\times I\times
M\to\Delta^p\times M$ by
$$
T^\lambda_{(m,\pi)}(t,0,x)=\tau^\lambda_{(m,\pi)}(t,x)
\quad\text{ and }\quad
T^\lambda_{(m,\pi)}(t,1,x)=\tau^{\lambda_1}_{(m,\pi)}(t,x)
$$
and extend it affine, where 
$\lambda_1=(\frac{1}{N},\dotsc,\frac{1}{N})\in C^\infty(M,\Delta^{N-1})$.
For $p\leq n$ we define a homotopy
$H_p:C_p^{{\Cal E}^\lambda_n,{\Cal U}}(B\bar G;\Z)\to 
C_{p+1}^{\Cal U}(B\bar G;\Z)$ on a $p$-simplex $\sigma$ by
$$
H_p(\sigma):=\sum_{i=1}^{p+1}(-1)^i\sum_{(m,\pi)\in
A^p_N}\sgn(\pi)(s_i^{p+1}\times\operatorname{id}_M)^*(T_{(m,\pi)}^\lambda)^*\sigma
$$

\proclaim
{Lemma~5.6}
In this situation ($n=2$) we have
$\partial\circ\varphi^\lambda_2=\varphi^\lambda_1\circ\partial$ and 
$(\varphi_1^\lambda)_*=j_*:H_1^{\Cal E^\lambda_2,\Cal
U}(B\bar G;\Z)\to H_1^{\Cal U}(B\bar G;\Z)$, where $j$ is the inclusion.
\endproclaim

{\it Proof}.
In this low dimension it is easy to see that 
$\partial\circ\varphi^\lambda_2=\varphi^\lambda_1\circ\partial$ and 
$\varphi_1^\lambda-\varphi_1^{\lambda_1}=\partial
H_1+H_0\partial$ (make a drawing!). 
So $(\varphi_1^\lambda)_*=(\varphi_1^{\lambda_1})_*$, but
the latter is the ordinary subdivison of the interval into $N$ subintervalls
and hence homotopic to the identity. So we obtain
$(\varphi_1^\lambda)_*=(\varphi_1^{\lambda_1})_*=j_*$.
\koniec

\remark
{Remark 5.7}
With some combinatorical difficulties one can show that this remains true
for arbitrary $n$.
\endremark

\remark
{Remark 5.8}
If $\lambda$ is subordinate to
an open cover $\Cal U$, and ${\Cal U}^{(n)}:=\{U_1\cup\cdots\cup
U_n:U_i\in{\Cal U}\}$ one easily sees
$$
\varphi^\lambda_p:C^{{\Cal E}^\lambda_n}_p(B\bar G;\Z)\to C^{{\Cal
U}^{(n)}}_p(B\bar G;\Z)
\qquad\forall p\leq n
$$
and therefore the name fragmentation mapping.
\endremark

\remark
{Remark 5.9}
If $\Cal V$ is an open covering with the property
$$
V_1,V_2\in\Cal V, V_1\cap V_2\neq\emptyset
\Rightarrow V_1\cup V_2\subseteq U\in\Cal U
$$
and $\lambda$ is subordinate to $\Cal V$ then $\varphi^\lambda_1$ induces a
mapping $(\varphi^\lambda_1)_*:H_1^{\Cal E^\lambda_2}(B\bar G;\Z)
\to H_1^{\Cal U}(B\bar G;\Z)$. One can see this as follows. If $\partial
d=c\in C_1^{\Cal E_2^\lambda}(B\bar G;\Z)$ with $d\in C_2^{\Cal
E_2^\lambda}(B\bar G;\Z)$ then Lemma~5.6 gives
$\varphi_1^\lambda(c)=\partial\varphi_2^\lambda(d)$, but we only have
$\varphi_2^\lambda(d)\in C^{\Cal V^{(2)}}_2(B\bar G;\Z)$. If we write
$$
\textstyle
\varphi_2^\lambda(d)=\sum_k\rho_k+\sum_l\kappa_l
$$
with $\su(\rho_k)\subseteq V_i\cup V_j$ for some $i,j$ with $V_i\cap
V_j\neq\emptyset$ and $\kappa_l$ such that there do not exist $V_i,V_j$ with
this property then one can show that $\partial(\sum_l\kappa_l)=0$ and hence
$\varphi_2^\lambda(d)=\partial(\sum_k\rho_k)$ with $\rho_k\in S^{\Cal
U}_2(B\bar G)$ by the construction of $\Cal V$.
\endremark

\proclaim
{Theorem 5.10}
Let $G$ be modular, $\Cal U$ be an open covering of $M$. Then the inclusion
induces an isomorphism $i_*:H^{\Cal U}_1(B\bar G;\Z)\cong H_1(B\bar G;\Z)$.
\endproclaim

{\it Proof}.
It suffices to show this for simplices which have support in a fixed compact
set $K$. Choose a covering $\Cal V$ with the property
$$
V_1,V_2\in\Cal V, V_1\cap V_2\neq\emptyset
\Rightarrow V_1\cup V_2\subseteq U\in\Cal U
$$
and choose $V_1,\dotsc, V_{N-1}\in\Cal V$ that cover $K$.
Let $\lambda\in C^\infty(M,\Delta^{N-1})$ be subordinated to 
$\{M\setminus K,V_1,\dotsc, V_{N-1}\}$ and let $r_*:H_1(B\bar G;\Z)\to
H_1^{\Cal E^\lambda_2}(B\bar G;\Z)$ be the inverse of the inclusion $j_*$
from Lemma~5.1. Then
$$
(\varphi_1^\lambda)_*\circ r_*:H_1(B\bar G;\Z)\to H_1^{\Cal
E^\lambda_2}(B\bar G;\Z)\to H_1^{\Cal U}(B\bar G;\Z)
$$
is an inverse of $i_*:H_1^{\Cal U}(B\bar G;\Z)\to H_1(B\bar G;\Z)$. Indeed
we have
$$
j_*=i_*\circ(\varphi_1^\lambda)_*:H_1^{\Cal E^\lambda_2}(B\bar G;\Z)
\to H_1^{\Cal U}(B\bar G;\Z)
\to H_1(B\bar G;\Z)
$$
by Lemma~5.6 with $\Cal U=\{M\}$ and so $i_*\circ(\varphi_1^\lambda)_*\circ
r_*=j_*\circ r_*=\operatorname{id}$. On the other hand the two mappings
$$
(\varphi_1^\lambda)_*\circ r_*\circ i_*:
H_1^{\Cal U}(B\bar G;\Z)
\to H_1(B\bar G;\Z)
\to H_1^{\Cal E_2^\lambda}(B\bar G;\Z)
\to H_1^{\Cal U}(B\bar G;\Z)
$$
and
$$
(\varphi_1^\lambda)_*\circ r_*:
H_1^{\Cal U}(B\bar G;\Z)
\to H_1^{\Cal U,\Cal E^\lambda_2}(B\bar G;\Z)
\to H_1^{\Cal U}(B\bar G;\Z)
$$
coincide and again by Lemma~5.6 $(\varphi_1^\lambda)_*\circ r_*\circ
i_*=(\varphi_1^\lambda)_*\circ r_*=j_*\circ r_*=\operatorname{id}$.
\koniec

\remark
{Remark 5.11}
With some combinatorical difficulties one can even show that the inclusion
induces isomorphisms $i_*:H^{\Cal U^{(n)}}_p(B\bar G;\Z)\cong H_p(B\bar
G;\Z)$ for all $p\leq n$, where $\Cal U^{(n)}:=\{U_1\cup\cdots\cup
U_n:U_i\in\Cal U\}$.
\endremark

\head
6. Proof of Theorem 1.1
\endhead

We begin with some consequences of Theorem 3.2.

\proclaim
{Proposition 6.1} Let $h_t$ be an isotopy in $\d$. Then $h_t$ can be
written as $h_t=h^s_t\cdots h^1_t$ where
$$
h^i_t=R_{\b ^i_t}R^{-1}_{\a}(g^i_t)^{-1}R_{\a}g^i_t\quad \forall t
$$
for some $\a \in T^k$ and some isotopies $g^i_t\in \d$ and $\b ^i_t\in T^k$.
\endproclaim

{\it Proof}. Observe first that $h_t$ can be written as the product of
$$h_{(p/m)t}h^{-1}_{(p-1/m)t},\ p=1,\ldots ,m,$$
 for $m$ sufficiently large.

Then we may assume that $h_t\in \Cal V$, where $\Cal V$ is a neighborhood of $id$
such that $R_{\alpha}\Cal V\subset \Cal U$, where $\Cal U, \alpha$ are as in
Theorem 3.2, and $\alpha $ is so small that $R_{\alpha }$ is in a
contractible neighborhood of $id$.
Thanks to Theorem 3.2 we have
 $R_{\alpha}h_t=R_{\b _t}g_t^{-1}R_{\alpha }g_t$ as required.
\Koniec

\proclaim
{Proposition 6.2} $\{R_{\b _t}\}=0$ in $H_1(\wt {\d})$.
\endproclaim

Indeed, it follows from Theorem 4 [10] and (3.1).

As a corollary we have
\proclaim
{Theorem 6.3}  $H_1(\wt {\d})=0$, that is $\wt{\d}$ and,
consequently, $\d$ are perfect.
\endproclaim

For the general case we need a more refined version of Theorem 6.3 which
can be formulated as follows.
\proclaim
{Proposition 6.4}
Let $U_2, V_2$ be  open balls in $T^{n-k}$ such that $\bar U_2\subset V_2$
 and let $V=T^k\times V_2$, $U=T^k\times U_2$.
 Suppose $\f|T^n-U$ is trivial, i.e. $L_y=T^k\t \{pt\}$ if $y\not \in U$.
 If $h_t$ is an isotopy
in $\du $ then $\{h_t\}=0$ in $H_1(\widetilde {\dv })$.
In other words, the map
$$
\i _*:H_1(\wt{\du})\r H_1(\wt {\dv})
$$
is trivial, where $\i :\du \r \dv$ is the canonical inclusion.
\endproclaim
{\it Proof}. Let $h_t\in \du$. In view of Theorem 6.3 we have
$$
h_t\sim [g^1_t,k^1_t]\ldots [g^r_t,k^r_t]$$
where $\sim $ denotes  the homotopy rel. endpoints, and $g_t^i,k_t^i\in
\d$.
We choose a smooth bump function $\mu :T^{n-k}\r [0,1]$ such that $\su \mu \subset V_2$
and $\mu =1$ on $U_2$. Then we let
$$ \bar g_t^i(x,y)=g^i_{\mu (y)t}(x,y),$$
$$ \bar k^i_t(x,y)= k^i_{\mu (y)t}(x,y),$$
where $(x,y)=(x_1,\ldots ,x_k,y_1,\ldots ,y_{n-k})$ is the standard chart
for $(T^n,\fk)$.
Observe that since $g^i_t$, $k_t^i$ are leaf preserving diffeomorphisms
then so are $\bar g^i_t$, $\bar k^i_t$.
 This follows from the assumption: $\f|V-U=\fk |V-U$.
Note as well that we
have
$$ h_t\sim [\bar g^1_t,\bar k^1_t]\ldots [\bar g^r_t,\bar k^r_t].$$
This is a consequence of the fact that
the initial homotopy is leafwise so that it can be modified in an obvious
manner.
            \Koniec

Now for technical reasons
we introduce the following sequence of open sets in $T^n$:
$$
\eqalign{
U&=U_1\times U_2\cr
U'&=U_1\times W_2\cr
V&=T^k\times V_2\cr
W&=T^k\times W_2\cr
W'&=T^k\times W'_2,\cr
}$$
where $U_1$ is an open ball in $T^k$,
and $U_2, V_2, W_2, W'_2$ are open balls in $T^{n-k}$ satisfying
$\bar U_2\subset V_2\subset \bar V_2\subset W_2 \subset \bar W_2\subset W'_2$.
 We have  the
following commutative diagram
$$
\CD
\wt{\du} @>\ti\i_1>> \wt {\dz}\\
@V\ti\i_3VV   @V\ti\i_2VV      \\
\wt{\dv} @>\ti\i_4>> \wt {\dw},\\
\endCD
$$
where  $\i_1:\du \r \dz$ is the canonical inclusion and so on.
The commutativity follows by the definition of the universal covering
and by the fact that $\i _j$ are inclusions. This diagram descends to
$$
\CD
H_1(\wt{\du}) @>\i_{1*}>>
H_1(\wt {\dz})\\
@V\iota _{3*}VV   @V\i _{2*}VV \\
H_1(\wt {\dv})  @>\iota _{4*}>>
H_1(\wt {\dw}).
\endCD
$$

Now let $\f'$ be a foliation on $M$
and let $h_t$ be an isotopy in
 $\operatorname {Diff}^{\infty}(M,\Cal F')_0$. By Proposition 2.4
 one can assume
that $h_t\in \du$ where $\f$ is some foliation on $T^n$ such that $\f|U=\f'|U$.
In addition, we assume that $(U',\phi)$ is a chart at $x$ chosen as in Theorem 2.2
and such that $U,U'$ are as above. This means that
$\phi (x)=0$, $dim(L_x)=k$, $\phi (U')=U_1\t W_2$, and  $\phi (\f'|U')=U_1\t \f'_2$
where $\f'_2$ is a foliation on $W_2$ with $L_x=\{pt\}$.

 The foliation $\f$ on $T^n$ can be defined as follows. $\f'_2$ being
 a smooth foliation on $W_2$, one has that its tangent distribution $T\f'_2$
 is determined by a family of smooth vector fields, say $S$ (cf.[13, p.545]).
 Then every vector field of $S$ respects $T\f'_2$, and $S$ is integrable
 [13, Cor.1]. Choose a bump function $\mu :T^{n-k}\r [0,1]$ with $\mu =1$
 on $U_2$ and $\su \mu \subset V_2$, where $V_2$ is as above. The
 family of vector fields $\mu S$ is still smooth and by the above
 reasoning it determines a foliation $\f_2$ on $T^{n-k}$. We let $\f=T^k\t \f_2$.
 Note that $\f=\fk$ in a neighborhood of $T^n-V$, and $\f=\f'$ on $U$.

Observe that
 Proposition 6.4 applies to $\f$ and we have that $\{h_t\}=0$
 in $H_1(\wt{\dv})$. Our
purpose is to show that $\{h_t\}=0$ in
the group $H_1(\wt {\dz})$. This will
be a consequence of the above diagram and Proposition 6.8 below.

 First we need the following two lemmata.

\proclaim
{Lemma 6.5}
With the above notation, there exist a finite family of open balls
 $\{W^i\}_{i=1}^s$ such that $W= \bigcup W^i$
and a related family of isotopies $\{\phi ^i_t\}_{i=1}^s$ in $Diff^{\inf}_{W'}(T^n,\f)_0$
such that $\phi ^i_1(W^i)\subset U'$
and
$$    \phi ^i_1|W^i\cap W^j\ =\ \phi_1^{ij}\circ \phi ^j_1|W^i\cap W^j\ ,$$
where $\phi _t^{ij}$ is an isotopy in $\dz$, for each $(i,j)$ such that $W^i\cap W^j\neq
 \emptyset$.       Moreover, we may have that $W^i\cap W^j$, whenever
 nonempty, is a ball.

\endproclaim
{\it Proof}.
There exists a covering $\{U^i_1\}_{i=1}^s$
of $T^k$ by open balls such that $U^i_1\cap U^j_1$ is a ball whenever
nonempty. Further, there exist  isotopies $\psi ^i_t$ in
$Diff^{\inf}(T^k)_0$, and $\psi ^{ij}_t$ in  $Diff^{\inf}_{U_1}(T^k)$
such that
$$
\psi ^i_t|U^i_1\cap U^j_1=\psi ^{ij}_1\circ \psi _1^j|U^i_1\cap U^j_1$$
(see, e.g.,{1]).
Then we let $W^i=U^i_1\times W_2$ and $\phi ^i_t=\psi ^i_{\mu (y)t}\times id$,
$\phi ^{ij}_t=\psi ^{ij}_{\mu (y)t}\times id$, where $\mu :T^{n-k}\r [0,1]$
is a bump function such that $\su \mu \subset W'_2$ and $\mu =1$ on $W_2$.
\koniec

Let us recall that $c\in B_n\bar G$ has its support in $U$ if and only if $\forall\
x,y\in \Delta ^n$ the diffeomorphism $c(x)c(y)^{-1}$ is supported in $U$.

The fragmentation property (Theorem 5.10) can be specified to our situation
as follows.

\proclaim
{Lemma 6.6} Let $U$, $W$, $\{W^i=U^i_1\t W_2\}_{i=1}^r$ be as above.
If a 1-chain $\sigma \in B \overline {\du}$ is a boundary of a 2-chain $c=\sum c_{\a} \in B\overline {\dv}$
then $\sigma $ is a boundary of a 2-chain $C=\sum C_{\a}$ such that the supports
of $C_{\a}$ are subordinate to $\{W^i\}$.
\endproclaim

For any $g\in G$, $G$ being a topological group, we denote by $I_g$ the
inner automorphism of $G$ induced by $g$.  Then
we have (cf.[2] or [1])
\proclaim
{Lemma 6.7} If $\sigma$, $\tau$ are any 1-simplices in $G$, then $\sigma$,
$I_{\tau(1)}\sigma$ are homological.
\endproclaim

Now we are in a position to prove

\proclaim
{Proposition 6.8} If $\sigma \in B_1\overline {\du }$ satisfies $\i _{2*}\i _{1*}\{\sigma \} =0$
then $\i _{1*}\{\sigma \}=0$.
\endproclaim

{\it Proof}. Let $\sigma \in B_1\overline {\du}$. Due to Proposition 6.4
 $\sigma =\p c$, where
$c=\sum c_{\a}\in B_2\overline {\dv}$.
Next, in light of Lemma 6.6 one has that support
of each $c_{\a}$ is contained in some $W^i$.

Under the  notation of Lemma 6.5 we assume
the convention :

(i) $\su (\partial _jc_{\a})\subset W^{i(j,\a )}$  and by $\phi_t^{i(j,\a )}$ we denote
the corresponding isotopy;

(ii) we assume that $W^{i(j,\a )}=U'$ and $\phi _t^{i(j,\a )}=id$, if $\su(\partial _jc_{\a})
\subset U'$;

(iii) if $\partial _jc_{\a}=\pm \partial _lc_{\beta }$ then $W^{i(j,\a)}=W^{i(l,\beta )}$ and
$\phi_t^{i(j,\a)}=\phi_t^{i(l,\beta )}$;

(iv) $\su(c_{\a })\subset W^{i(\a)}$ and $\phi_t^{i(\a)}$ denotes the corresponding
isotopy.

We have the following equality:
$$ \sigma \ =\ \sum _{\a }\sum _{j=0}^2\ (-1)^j\partial _jc_{\a}\ =\
\sum _{\a }\sum _{j=0}^2\ (-1)^jI_{\phi _1^{i(j,\a )}}(\partial _jc_{\a }) .\leqno (6.1)$$
If fact, if support of the edge $\partial _jc_{\a }$ is in $U'$ then,
because of (ii), nothing changes in the r.h.s.
Otherwise, this edge must be reduced in the sum on the l.h.s., and by
(iii) so must be on the r.h.s.

Due to Lemmata 6.5 and 6.7 we get
$$
I_{\phi _1^{i(j,\a)}}(\partial _jc_{\a})\ =\ I_{\phi ^{i(j,\a)i(\a)}}I_{\phi_1^{i(\a )}}(\partial _jc_{\a })\ \sim\
I_{\phi_1^{i(\a )}}(\partial _jc_{\a })\ =\ \partial _jI_{\phi_1^{i(\a )}}(c_{\a }),$$
where $\sim$ stands for the homology relation.
Combining this with (6.1) we have
$$  \sigma\ \sim \ \sum_{\a }\sum_{j=0}^2\ (-1)^j\partial _jI_{\phi _1^{i(\a )}}(c_{\a })\ =\
\partial \sum_\a I_{\phi _1^{i(\a)}}(c_{\a }) .$$
In view of (iv), $\su(\sum _{\a} \ I_{\phi _1^{i(\a )}}(c_{\a}))\subset U'$.
Thus $\sigma $ is a coboundary in $B\overline {\dz}$, i.e. $\i_{1*}\{\sigma \}
=0$.
\Koniec

By combining Propositions 6.4 and 6.8 it is visible  that Theorem 1.1
holds.

\head
References
\endhead

[1] A.Banyaga, The Structure of Classical Diffeomorphism Groups,
Kluwer Academic Publishers, Dordrecht 1997.

[2] K.S.Brown, Cohomology of Groups, Springer Verlag, New York - Heidelberg -
Berlin 1982.

[3] P.Dazord, {\it Feuilletages \`a singularit\'es}, Indag. Math. {\bf 47}
(1985), 21-39.

[4] D.B.A.Epstein, {\it The simplicity of certain groups of homeomorphisms},
Comp. Math. {\bf 22}(1970), 165-173.

[5] K.Fukui, {\it Homologies of the group $Diff^{\infty }(R^n,0)$ and its
subgroups}, J. Math. Kyoto Univ. {\bf 20}(1980), 475-487.

[6] M.R.Herman, {\it Sur le groupe des diff\'eomorphismes du tore}, Ann.
Inst. Fourier (Grenoble) {\bf 23}(1973), 75-86.

[7] A.Kriegl, P.W.Michor, {\it The Convenient Setting of Global Analysis},
American Mathematical Society, 1997.

[8]  J.N.Mather, {\it A curious remark concerning the geometric transfer
map}, Comment. Math. Helv. {\bf 59}(1984), 86-110.

[9]  J.P.May, Simplicial objects in algebraic topology, Van Nostrand 1967.

[10] T.Rybicki, {\it The identity component of the leaf preserving diffeomorphism
group is perfect}, Mh. Math. {\bf 120}(1995), 289-305.

[11] T.Rybicki, {\it On the group of Poisson diffeomorphisms of the
 torus}, Rendiconti Mat. (Roma), in press.

[12] P.Stefan, {\it Accesibility and foliations with singularities}, Bull.
Amer. Math. Soc. {\bf 80}(1974), 1142-1145.

[13] P.Stefan, {\it Integrability of systems of vectorfields}, J. London
Math. Soc. {\bf 21}(1980), 544-556.

[14] F.Sergeraert, {\it Un th\'eor\`eme de fonctions implicites sur
certains espaces de Fr\'echet et quelques applications}, Ann. Scient. Ec.
Norm. Sup. 4 s\'erie  {\bf 5}(1972), 599-660.

[15] W.Thurston, {\it Foliations and groups of diffeomorphisms,} Bull. Amer.
Math. Soc. {\bf 80}(1974), 304-307.

\vskip1truecm
(T.R.)

Institute of Mathematics

Pedagogical University

ul.Rejtana 16 A

35-310 Rzesz\'ow, POLAND

e-mail: rybicki\@im.uj.edu.pl

\vskip1truecm
(St.H.)

Institute of Mathematics

University of Vienna

Strudlhofgasse 4

1090 Wien, Austria

e-mail: stefan\@nelly.mat.univie.ac.at
\enddocument